%% file: GeneralCumulant2.tex
\documentclass[12pt,twoside]{amsart}

\input{macros.tex}

\newcommand{\externalref}[2]{#1.\ref{I-#2}}
\newcommand{\ssection}[1]{\section{#1}}

\begin{document}

\title[Cumulants in Noncommutative Probability Theory II]
{Cumulants in Noncommutative Probability Theory II. Generalized Gaussian random variables}
\author{Franz Lehner}

\thanks{Supported by the European Network \no{}HPRN-CT-2000-00116
and the Austrian Science Fund (FWF) Project \no{}R2-MAT}

\address{
Franz Lehner\\
In\-sti\-tut f\"ur Mathe\-ma\-tik C\\
Tech\-ni\-sche Uni\-ver\-si\-t\"at Graz\\
Stey\-rer\-gas\-se 30, A-8010 Graz\\
Austria}
\email{lehner@finanz.math.tu-graz.ac.at}


\date{\today}

\begin{abstract}
  We continue the investigation of noncommutative cumulants.
  In this paper various characterizations of generalized Gaussian random
  variables are proved.
\end{abstract}

\maketitle{}
\ddate{27.11.2000}

\tableofcontents{}

Generalized Gaussian random variables and Brownian motions have a long
history in noncommutative probability theory and noncommutative
central limit theorems.
For a systematic study see~\cite{GutaMaassen:2002:generalised}.
In this paper we consider generalized Gaussian random variables
from the point of view of combinatorial cumulant theory
as developed in our paper~\cite{Lehner:2002:Cumulants1},
to which we refer as part~I.
Our aim is to prove characterizations of Gaussian random variables,
as found in
\cite{KaganLinnikRao:1973:characterization} and 
\cite{Bryc:1995:normal} for classical Gaussian distributions.
There are essentially two kinds of characterizations.
The proofs of the simpler ones like spherical symmetry, Bernstein's and
Lukacs' theorems can be immediately transferred to the noncommutative case,
while other theorems, notably including Cram\'er's and Marcinkiewicz' theorems,
do not hold in general.

\input{GeneralCumulantGaussian.tex}

\bibliography{GeneralCumulants}
\bibliographystyle{mamsalpha}
\end{document}

%% file: macros.tex
\usepackage[T1]{fontenc}
\usepackage[latin1]{inputenc}
\usepackage{amsmath,amsfonts,amssymb,amscd,amsthm,amsxtra}

\usepackage{a4wide}
\usepackage{url}
\usepackage{enumerate,verbatim}
\usepackage{ifthen}

\usepackage{url,textcomp}
\newcommand{\no}{\textnumero}

\usepackage{xr1}

\usepackage{pstricks}
\newcommand{\ddate}[1]{}

\newcommand{\ppaper}[4][ ]{}
\newcommand{\paper}[4][ ]{}

\numberwithin{equation}{section}

\makeatletter
  \DeclareRobustCommand\em
         {\@nomath\em \ifdim \fontdimen\@ne\font >\z@
                       \upshape \else \slshape \fi}
\makeatother

\newtheoremstyle{mythm}
  {9pt}
  {9pt}
  {\slshape}
  {0pt}
  {\bfseries}
  {.}
  { }
  {\thmname{#1} \thmnumber{#2}\thmnote{ (#3)}}

\theoremstyle{mythm} 
\newtheorem{Theorem}{Theorem}[section]
\newtheorem{Lemma}[Theorem]{Lemma}
\newtheorem{Proposition}[Theorem]{Proposition}

\newtheorem{Corollary}[Theorem]{Corollary}

\theoremstyle{definition} 
\newtheorem{Definition}[Theorem]{Definition}

\newtheorem{Remark}[Theorem]{Remark}

\newcommand{\file}[1]{{\url{#1}}}

\newcommand{\IC}{\mathbf{C}}                     
\newcommand{\IN}{\mathbf{N}}                     
\newcommand{\IR}{\mathbf{R}}                     
\newcommand{\SG}{\mathfrak{S}}                   
\newcommand{\eps}{{\varepsilon}}                 
\renewcommand{\phi}{{\varphi}}                   
\newcommand{\abs}[1]{\left\lvert #1 \right\rvert}  
\newcommand{\exchm}{\alg{E}} 
\newcommand{\exchF}{{\alg{F}}} 
\newcommand{\exch}{\mbox{$\exchm$}} 

\newlength{\tmpl}

\newcommand{\indep}{\perp\kern-3pt\perp}

\newcommand{\alg}[1]{{\mathcal{#1}}}                
\newcommand{\OrthogonalGroup}{\mathcal{O}}

\DeclareMathOperator{\tr}{{\mathrm tr}}                     



%% file: GeneralCumulantGaussian.tex

\ssection{Preliminaries}
Throughout this paper we will consider a fixed noncommutative probability space~$(\alg{A},\phi)$ 
and an exchangeability system~$\exchm= (\alg{U},\tilde{\phi},\alg{J})$
for~$(\alg{A},\phi)$
as defined in part~I.
The interchangeable images of~$\alg{A}$ in $\alg{U}$ will as usual be denoted 
by~$(\alg{A}_i)_{i\geq0}$, and we shall identify~$\alg{A}$ with $\alg{A}_0$.
We shall moreover assume that $\alg{A}$ is a $*$-algebra and
that all considered random variables are selfadjoint.
\begin{Definition}
  \label{def:GeneralCumulants:Gaussian:SameDistribution}
  We say that two random variables~$X$ and~$Y\in\alg{A}$
  \emph{have the same distribution given~$\exchm$},
  if for any word $W=W_1W_2\cdots W_n$ with
  $W_i\in \{X\}\cup \bigcup_{i\geq1} \alg{A}_i$ the expectation~$\tilde{\phi}(W)$
  does not change if we replace each occurrence of~$X$ by~$Y$.
  We call~$X$ and~$Y$ \emph{\exch-i.i.d.}\, if in addition they are \exch-independent.
  Similarly a sequence $(X_i)_{i\in\IN}\subseteq\alg{A}$ of \exch-independent random variables
  is called \emph{\exch-i.i.d.}\
  if for any word $W=W_1W_2\cdots W_n$ with
  $W_i\in \{X_i: i\in\IN\}\cup \bigcup_{i\geq1} \alg{A}_i$ the expectation~$\tilde{\phi}(W)$
  does not change if we apply a permutation $\sigma\in\SG_\infty$ to the indices
  of $X_i$, i.e., if we replace each occurrence of $X_i$ by $X_{\sigma(i)}$.
\end{Definition}
We will need the following weak variant of pyramidal independence
(cf.~Definition~\externalref{I}{def:GeneralCumulants:PyramidalIndependence}).
\begin{Definition}
  \label{def:GeneralCumulants:Gaussian:SingletonCondition}
  Let $X_i$ be an \emph{interchangeable sequence} of (centered) random variables,
  that is, for every permutation $\pi\in\SG_\infty$ and every choice of indices
  $i_1,i_2,\dots,i_n$ the expectation does not change under permutations:
  $$
  \phi(X_{\pi(i_1)} X_{\pi(i_2)} \cdots X_{\pi(i_n)}) = 
  \phi(X_{i_1} X_{i_2} \cdots X_{i_n})
  .
  $$
  We say that the \emph{singleton condition} holds if
  $$
  \phi(X_{i_1} X_{i_2} \cdots X_{i_n}) = 0
  $$
  whenever one of the $X_i$'s occurs exactly once.
\end{Definition}

\ddate{24.01.2001}

Let us start this section by quoting a general noncommutative central limit
theorem. 

\begin{Theorem}[{\cite{BozejkoSpeicher:1996:interpolations}}]
  \label{thm:GeneralCumulants:CentralLimitTheorem}
  Let $(\alg{A}, \phi)$ be a noncommutative probability space,
  and $X_i=X_i^*\in \alg{A}$ be a sequence of exchangeable random variables.
  For a partition $\nu$ denote
  $$
  \phi(\nu) = \phi(X_{i_1}X_{i_2}\dotsm X_{i_n})
  $$
  where $(i_1,i_2,\dots,i_n)$ is any multiindex with kernel $\nu$.
  Assume that $\phi(X_i)=0$, $\phi(X_i^2)=1$ and moreover that
  the \emph{singleton condition} of
  Definition~\ref{def:GeneralCumulants:Gaussian:SingletonCondition} holds.
  Then the sequence $S_N = \frac{1}{\sqrt{N}} \sum_1^N X_i$
  has limit distribution
  $$
  \lim_{N\to \infty} \phi(S_N^{2n+1}) = 0 
  \qquad  \qquad
  \lim_{N\to \infty} \phi(S_N^{2n})
  = \sum_{\nu\in\Pi_{2n}^{pair}} \phi(\nu)
  $$
\end{Theorem}

Interchangeable sequences generate interchangeable algebras and give rise to exchangeability
systems.
In view of the preceding noncommutative central limit theorem
we define Gaussian families as follows
(see also \cite{GutaMaassen:2002:generalised}).
\begin{Definition}
  \label{def:GeneralCumulants:Gaussian}
  An interchangeable family $(X_i)$ of random variables is called (centered) \emph{Gaussian}
  if all cumulants which correspond to non-pair partitions vanish.
  In other words, there is a function on pair partitions
  $\nu:\Pi_n^{(2)}\to\IC$ such that for all $h:[n]\to I$
  \begin{equation}
    \label{eq:GeneralCumulants:GaussDefinition}
    \phi(X_{h(1)} X_{h(2)} \dotsm X_{h(n)})
    = \sum_{\substack{\pi\in\Pi_n^{(2)}\\ \pi\leq\ker h }}
       \nu(\pi)
  \end{equation}
  In particular, odd moments vanish and the singleton condition holds.
\end{Definition}
Noncommutative (i.e.\ operator valued)
Khinchin inequalities are available for Gaussian families, see
\cite{Buchholz:2001:operator}.

In the following all random variables are assumed self-adjoint
and the state $\tilde{\phi}$ is assumed to be faithful. 
This is needed for the following crucial lemma to be valid.
\begin{Lemma}
  \label{lem:GeneralCumulantGaussian:phifaithfulimpliesK2positive}
  Let $\exchm=(\alg{U},\tilde{\phi},\alg{J})$ be an exchangeability system
  for a $C^*$-probability space~$(\alg{A},\phi)$
  with $\tilde{\phi}$ faithful.
  Let $X\in\alg{A}$ be such that $X^{(1)}\ne X^{(2)}$.
  Then $K^\exchm_2(X,X)>0$.
\end{Lemma}
\begin{proof}
  By Good's formula
  $$
  K^\exchm_2(X,X)
  = \phi((X^{(1)})^2 - X^{(1)}X^{(2)})
  = \frac{1}{2} \phi( (X^{(1)} - X^{(2)} )^2 ) > 0
  $$
\end{proof}

\ssection{Spherical symmetry and related characterizations}

We review now some characterizations of classical Gaussians which
may or may not hold in the general framework of definition~\ref{def:GeneralCumulants:Gaussian}.
A simple characterization is the following.
\begin{Proposition}
  \label{prop:GeneralCumulantGaussian:XsameasXpYsqrt2impliesgaussian}
  Let $X$ and $Y$ be \exch-i.i.d.\ 
  noncommutative random variables such that
  $\frac{X+Y}{\sqrt{2}}$ has the same distribution as $X$ (and $Y$).
  Then $X$ and $Y$ are centered Gaussian.
\end{Proposition}
\begin{proof}
  We have to show that
  for every partition $\pi$ the cumulant $K^\exchm_\pi(X,X,\dots,X)$ vanishes
  unless $\pi$ is a pair partition.
  Indeed, whenever there is a block~$B$ of length $m\ne2$,
  then the assumption and Good's formula
  (Proposition~\externalref{I}{prop:GeneralCumulants:PartitionGoodFormula}) imply
  $$
  K^\exchm_\pi(X,X,\dots,X) = K^\exchm_\pi(X_1,X_2,\dots,X_n)
  $$
  where 
  $$
  X_i=
  \begin{cases}
    \frac{X+Y}{\sqrt{2}} & \text{if $i\in B$}\\
    X                    & \text{if $i\not\in B$}
  \end{cases}
  $$
  and because of vanishing of mixed cumulants we obtain
  $$
  K^\exchm_\pi(X_1,X_2,\dots,X_n)
  = 2^{-m/2} K^\exchm_\pi(X,X,\dots,X) + 2^{-m/2}K^\exchm_\pi(X_1',X_2',\dots,X_n')
  $$
  where 
  $$
  X_i'=
  \begin{cases}
    Y & \text{if $i\in B$}\\
    X & \text{if $i\not\in B$}
  \end{cases}
  $$
  and thus
  $$
  K^\exchm_\pi(X,X,\dots,X) = 2^{1-m/2} K^\exchm_\pi(X,X,\dots,X)
  $$
  which is only possible if the cumulant vanishes.
\end{proof}

The following generalization holds. See~\cite[Thm.~3.3.1]{Bryc:1995:normal} for the
classical case.
\begin{Proposition}
  Let $X_i$ be a sequence of \exch-i.i.d. random variables and suppose that there
  are real numbers~$a_1,a_2,\dots,a_n\ne0$ such that $\sum a_i X_i$ has the same
  distribution as~$X=X_1$. Then $X$ is gaussian.
\end{Proposition}
\begin{proof}
  Comparing the second cumulant we get
  $$
  K^\exchm_2(X) = \sum a_i^2 K^\exchm_2(X)
  $$
  and by Lemma~\ref{lem:GeneralCumulantGaussian:phifaithfulimpliesK2positive}
  we infer that $\sum a_i^2=1$.
  The rest of the proof is the same as the proof of
  Proposition~\ref{prop:GeneralCumulantGaussian:XsameasXpYsqrt2impliesgaussian}.
\end{proof}

Another characterization of classical Gaussians is \emph{Maxwell's theorem}.
Its analogue for classical exchangeable random variables was proved by
\cite{Kingman:1972:random}%
\paper{Kingman}{On random sequences with spherical symmetry}{},
namely that spherically symmetric
exchangeable classical random variables are conditionally i.i.d.\ normal.
For the free case see \cite{Nica:1996:Rtransforms}%
\ppaper{Nica}{${R}$-transforms of free joint distributions and non-crossing partitions}{}.

\begin{Definition}
  A family $(X_i)$ of random variables is \emph{spherically symmetric}
  if for every $n\in\IN$ and for every real orthogonal matrix
  $U\in \OrthogonalGroup(n)$ the families
  $Y_i = \sum U_{ij} X_j$ and $X_i$ have the same joint distribution.
\end{Definition}

\begin{Theorem}
  \label{thm:GeneralCumulants:Maxwell}
  An infinite interchangeable family $(X_i)$ is Gaussian if and only if it is spherically
  symmetric.
\end{Theorem}

\begin{proof}
  Assume $X_i$ is gaussian and let $U = [U_{ij}]\in \OrthogonalGroup(n)$ be an
  arbitrary orthogonal matrix, i.e.,
  $\sum_j U_{ij} U_{kj} = \delta_{ik}$.
  Fix an index map $g:[m]\to [n]$.
  Then by multilinearity we have
  \begin{align*}
    \phi(Y_{g(1)} Y_{g(2)} \dotsm Y_{g(m)})
    &= \sum_{h:[m]\to[n]}
        U_{g(1),h(1)} U_{g(2),h(2)}\dotsm  U_{g(m),h(m)}
        \,
        \phi(X_{h(1)} X_{h(2)} \dotsm X_{h(m)}) \\
    &= \sum_{h:[m]\to[n]}
        U_{g(1),h(1)} U_{g(2),h(2)}\dotsm  U_{g(m),h(m)}
        \,
        \sum_{\substack{\pi\in\Pi_m^{(2)}\\ \pi\leq\ker h }}
         \nu(\pi)\\
\intertext{this is zero unless $m$ is even and in the latter case}
    &= \sum_{\pi\in\Pi_m^{(2)}}
        \nu(\pi)
        \sum_{\substack{h:[m]\to[n]\\ \ker h\geq \pi}}
        U_{g(1),h(1)} U_{g(2),h(2)}\dotsm  U_{g(m),h(m)}\\
    &= \sum_{\pi\in\Pi_m^{(2)}}
        \nu(\pi)
        \prod_{\{a,b\}\in \pi}
        \sum_j U_{g(a),j} U_{g(b),j}\\
    &= \sum_{\pi\in\Pi_m^{(2)}}
        \nu(\pi)
        \prod_{\{a,b\}\in \pi}
        \delta_{g(a),g(b)} \\
    &= \sum_{\pi\leq \ker g}
        \nu(\pi)\\
    &= \phi(X_{g(1)} X_{g(2)} \dotsm X_{g(m)})
  \end{align*}
  
  For the converse, we are going to prove that if
  the $X_i$ are even and if for each~$n$ there exists a 
  orthogonal matrix, none of whose entries has modulus~$1$
  and  which leaves the joint distribution invariant,
  then only the cumulants of order~$2$ are nonzero.

  Suppose $(X_i)_{i=1,2,\dots,n}$ has the same distribution as
  $Y_i = \sum  U_{ij} X_j$ 
  for some orthogonal matrix with the above property.
  We have to show that the cumulants corresponding to non-pair partitions
  vanish. Let $\pi=\{\pi_1,\pi_2,\dots,\pi_p\}\in\Pi_m$.
  Because the $X_i$ are even,
  we can easily dispose of partitions with singletons, so assume that all
  blocks have cardinality at least $2$.
  Then by multilinearity
  \begin{align*}
    K^\exchm_{\pi}(X,X,\dots,X)
    &= K^\exchm_{\pi}(X_{\pi(1)},X_{\pi(2)},\dots,X_{\pi(m)})\\
    &= K^\exchm_{\pi}(Y_{\pi(1)},Y_{\pi(2)},\dots,Y_{\pi(m)})\\
    &= \sum_{h:[n]\to[n]}
        U_{\pi(1),h(1)} U_{\pi(2),h(2)}\dotsm  U_{\pi(m),h(m)}
        K^\exchm_\pi(X_{h(1)}, X_{h(2)}, \dots, X_{h(m)})\\
    &= \sum_{\ker h\geq \pi}
        U_{\pi(1),h(1)} U_{\pi(2),h(2)} \dotsm U_{\pi(m),h(m)}
        K^\exchm_{\pi}(X,X,\dots,X)\\
    &= \prod_{\pi_j\in\pi}
        \left(
          \sum_k U_{j,k}^{\abs{\pi_j}}
        \right)
        K^\exchm_{\pi}(X,X,\dots,X)\\
  \end{align*}
  and by assumption each $\sum_k U_{jk}^{\abs{\pi_j}}$
  where $\abs{\pi_j}>2$ has modulus strictly less than $1$
  and thus the product is different from $1$ unless all block sizes are equal to $2$.
\end{proof}

Actually a stronger characterization holds,
known as \emph{Bernstein's theorem} or \emph{Kac-Lo\`eve theorem}
 \cite[\S{}III.4]{Feller:1971:introductionII}.
\begin{Proposition}
  Let~$X_1$ and $X_2$ be \exch-independent noncommutative random variables
  and assume that the random variables~$Y_1=\alpha X_1 + \beta X_2$
  and~$Y_1=\gamma X_1 + \delta X_2$ are also \exch-independent
  with $\alpha\gamma+\beta\delta=0$ and $\alpha,\beta,\gamma,\delta\ne0$
  (that is, the matrix
  $
  \left[
    \begin{smallmatrix}
      \alpha & \beta\\
      \gamma & \delta
    \end{smallmatrix}
  \right]
  $
  has orthogonal columns).
  Then $X_1$ and $X_2$ are (shifted) gaussian and have the same variance.
\end{Proposition}
\begin{proof}
  The \exch-independence of $Y_1$ and $Y_2$ implies vanishing of mixed cumulants,
  in particular
  \begin{align*}
    0&= K^\exchm_2(Y_1,Y_2)
      = \alpha \gamma K^\exchm_2(X_1,X_1)
        +\beta \delta K^\exchm_2(X_2,X_2) \\
     &= \alpha\gamma( K^\exchm_2(X_1,X_1) - K^\exchm_2(X_2,X_2) )
  \end{align*}
  and therefore $K^\exchm_2(X_1,X_1)=K^\exchm_2(X_2,X_2)$.
  Moreover we have
  \begin{align*}
    \begin{bmatrix}
      0\\
      0
    \end{bmatrix}
    &= \begin{bmatrix}
         K^\exchm_n(Y_1,Y_1,\dots,Y_1,Y_1,Y_2) \\
         K^\exchm_n(Y_1,Y_1,\dots,Y_1,Y_2,Y_2) 
       \end{bmatrix}\\
    &= \begin{bmatrix}
         \alpha^{n-1} \gamma K^\exchm_n(X_1,X_1,\dots,X_1)
           + \beta^{n-1} \delta K^\exchm_n(X_2,X_2,\dots,X_2) \\
         \alpha^{n-2} \gamma^2 K^\exchm_n(X_1,X_1,\dots,X_1)
           + \beta^{n-2} \delta^2 K^\exchm_n(X_2,X_2,\dots,X_2) 
       \end{bmatrix} \\
    &= \begin{bmatrix}
         \alpha & \beta \\
         \gamma & \delta
       \end{bmatrix}
       \begin{bmatrix}
         \alpha^{n-2}\gamma K^\exchm_n(X_1,X_1,\dots,X_1)\\
         \beta^{n-2} \delta K^\exchm_n(X_2,X_2,\dots,X_2) 
       \end{bmatrix}
  \end{align*}
  and since the matrix 
  $
  \left[
    \begin{smallmatrix}
      \alpha & \beta\\
      \gamma & \delta
    \end{smallmatrix}
  \right]
  $
  is invertible (if not, the random variables $Y_1$ and $Y_2$ are actually scalar multiples of
  each other and cannot be \exch-independent),
  the higher order cumulants vanish.
\end{proof}

For the multidimensional version of Bernstein's theorem we need the following
class of matrices.

\begin{Definition}
  A matrix $A\in M_n(\IR)$ is called \emph{reducible} if there are permutation
  matrices $C_1$,~$C_2$ such that
  $$
  C_1 A C_2 = 
  \begin{bmatrix}
    A_1 & 0\\
    0   & A_2
  \end{bmatrix}
  $$
  Otherwise~$A$ is called irreducible. Equivalently,
  $A$ is irreducible if it does not commute with any projection of the form
  $P=\sum_{i\in I} e_{ii}$.
\end{Definition}

\begin{Proposition}[{\cite[Thm.~3.5]{HiwatashiNagisaYoshida:1999:characterizations}}]
  Let $(\alg{A},\phi)$ be a $C^*$-probability space with $\phi$ faithful.
  Let~$X_i\in\alg{A}$, $i=1,2,\dots,n$ (with $n\geq3$)
  be centered \exch-independent random variables and
  let $U$ be an irreducible orthogonal $n\times n$ matrix
  such that $Y_i=\sum U_{ij}X_j$ are also \exch-independent.
  Then the~$X_i$ are \exch-i.i.d.\ Gaussian.
\end{Proposition}
\begin{proof}
  First let us prove that all~$X_i$ have the same variance.
  Indeed,
  \begin{align*}
    K^\exchm_2(Y_i,Y_j)
    &= \sum_{k,l} U_{ik} U_{jl} K^\exchm_2(X_k,X_l) \\
    &= \sum_k U_{ik} U_{jk} K^\exchm_2(X_k,X_k) \\
  \end{align*}
  i.e., if we set $\Xi = [ K^\exchm_2(X_i,X_j) ]$, $H=[ K^\exchm_2(Y_i,Y_j)]$
  (both are diagonal matrices by assumption and have nonzero diagonal entries
  by Lemma~\ref{lem:GeneralCumulantGaussian:phifaithfulimpliesK2positive}),
  then we have
  $$
  H = U \Xi U^t
  .
  $$
  Because the spectrum is invariant, it follows that we can also write $H$ as
  a permutation of $\Xi$: $H=C^t\Xi C$.
  Consequently $CU\Xi = \Xi CU$ and~$CU$ is irreducible and commutes with the spectral
  projections of $\Xi$. 
  The latter have the form~$\sum_{i\in I} e_{ii}$ and therefore~$\Xi$ is a multiple
  of the identity matrix.
  To conclude the proof we have to show that the higher order cumulants vanish.
  Note that every row of $U$ has at least two nonzero entries.
  (If there is only one nonzero entry, the other rows must have zero in the
  corresponding entry because of orthogonality, causing the matrix to be reducible).
  We fix an index~$k$ and assume without loss of generality that 
  the entries $U_{1,k}$ and $U_{2,k}$ are nonzero.
  Consider for $m\geq3$ the identity
  \begin{align*}
    0&= K^\exchm_m(Y_i,Y_1,Y_1,\dots,Y_1,Y_2) \\
     &= \sum_j U_{ij} U_{1j}^{m-2} U_{2j} K^\exchm_m(X_j,X_j,\dots,X_j) 
  \end{align*}
  which holds for every~$i$.
  Therefore we have by orthogonality
  \begin{align*}
    0&= \sum_i U_{ik}
         \sum_j
          U_{ij} U_{1j}^{m-2} U_{2j} K^\exchm_m(X_j,X_j,\dots,X_j) \\
    0&= \sum_j \delta_{jk}
         U_{1j}^{m-2} U_{2j} K^\exchm_m(X_j,X_j,\dots,X_j) \\
     &=  U_{1k}^{m-2} U_{2k} K^\exchm_m(X_k,X_k,\dots,X_k)
  \end{align*}
  It follows that $K^\exchm_m(X_k,X_k,\dots,X_k) = 0$.
\end{proof}

\ssection{Linear Forms. The Skitovi\v{c}-Darmois theorem and its relation to
  Cram\'er's  and Marcinkiewicz' theorem}

In the classical case, an even stronger result than Bernstein's theorem holds,
known as \emph{Skitovi\v{c}-Darmois theorem}.
\begin{Theorem}[{\cite[Ch.~3]{KaganLinnikRao:1973:characterization}}]
  Let for $n\ge2$ classical independent random variables
  $X_1,X_2,\dots,X_n$ be given
  and let $a_i$, $b_i$ be real numbers for which $a_ib_i\ne0$
  for each $i$. Assume that the linear statistics
  $$
  Y_1 = a_1 X_1 + a_2 X_2 + \dots + a_n X_n
  \qquad
  Y_2 = b_1 X_1 + b_2 X_2 + \dots + b_n X_n
  $$
  are independent. Then $X_i$ are all gaussian.
\end{Theorem}
This theorem heavily depends on Marcinkiewicz' and Cram\'er's theorems.
\begin{Theorem}[{Marcinkiewicz \cite{Marcinkiewicz:1939:propriete,Bryc:1995:normal}}]
  Let $X$ be a classical random variable with only finitely many non-vanishing 
  (classical) cumulants.
  Then $X$ is normal, i.e., all cumulants of order greater than $2$ vanish.
\end{Theorem}
\begin{Theorem}[{Cram\'er \cite[\S XV.8]{Feller:1971:introductionII}}]
  Let $X_1$,\dots, $X_n$ be classical independent random variables
  such that their sum $X_1+X_2+\dots+X_n$ is normal.
  Then all $X_i$ are normal.
\end{Theorem}
In the general case Marcinkiewicz' and Cram\'er's theorems do not hold,
for example in free probability \cite{BercoviciVoiculescu:1995:superconvergence}.
Counterexamples to both can be fabricated from the following theorem,
which shows that the free cumulants of order higher than $2$ can take more
or less arbitrary values.
\begin{Theorem}[{\cite[Thm.~2]{BercoviciVoiculescu:1995:superconvergence}}]
  For every $r>0$ there exists $\delta>0$
  such that the Taylor coefficients $c_n$
  of every function $f(z)=-z+\sum_{n=0}^\infty c_{n+1} z^n$
  analytic in $\{z:\abs{z}<r\}$ which satisfies
  $f(\bar z)=\overline{f(z)}$ and $\abs{f(z)}<\delta$ for every $z$
  are the free cumulants of a probability measure.
\end{Theorem}
\begin{Corollary}
\label{cor:GeneralCumulants:FreeCramerFailure}
  For small enough $\eps$ there exists a selfadjoint random variable $X$
  (equivalently, a probability measure on the real line)
  with free cumulants $K^\exchF_1(X)=0$, $K^\exchF_2(X)=1$, $K^\exchF_3(X)=\eps$ 
  and $K^\exchF_n(X)=0$ for $n\geq4$.
\end{Corollary}
This lemma can be used to show that the analogue of Skitovi\v{c}' theorem
fails in the free case if there are at least three random variables involved.
\begin{Proposition}
  There are free random variables $X_1,X_2,X_3$ which are not semicircular
  and such that $Y_1=a_1X_1+a_2X_2+a_3X_3$ and $Y_2=b_1X_1+b_2X_2+b_3X_3$ 
  are free.
\end{Proposition}
\begin{proof}
  By Corollary~\ref{cor:GeneralCumulants:FreeCramerFailure} there exists
  $\eps_0>0$ such that for every $\eps$ with $\abs{\eps}<\eps_0$ there exists
  a selfadjoint random variable $X(\eps)$ such that all free cumulants of 
  $X(\eps)$ are zero with the exceptions $K^\exchF_2(X(\eps))=1$ and
  $K^\exchF_3(X(\eps))=\eps$.
  Let $X_1,X_2,X_3$ be a free family where
  $X_1\sim X(\eps/4)$,   $X_2\sim X(\eps)$ 
  and $X_3\sim X(\eps)$, 
  $$
  Y_1 = 2 X_1 - X_2 + 2 X_3
  \qquad
  Y_2 = 2 X_1 +2 X_2 - X_3
  $$
  Then $Y_i$ have vanishing mixed cumulants:
  \begin{align*}
    K^\exchF_2(Y_1,Y_2) &= a_1 b_1 K^\exchF_2(X_1,X_1)
                    + a_2 b_2 K^\exchF_2(X_2,X_2) + a_3 b_3 K^\exchF_2(X_3,X_3) \\
                 &= 4 - 2 - 2 = 0 \\
    K^\exchF_3(Y_1,Y_1,Y_2) &= a_1^2 b_1 K^\exchF_3(X_1,X_1,X_1)
                        + a_2^2 b_2 K^\exchF_3(X_2,X_2,X_2)
                        + a_3^2 b_3 K^\exchF_3(X_3,X_3,X_3) \\
                     &= 8\cdot\frac{1}{4} + 2  - 4 = 0 \\
    K^\exchF_3(Y_1,Y_2,Y_2) &= a_1 b_1^2 K^\exchF_3(X_1,X_1,X_1)
                        + a_2 b_2^2 K^\exchF_3(X_2,X_2,X_2)
                        + a_3 b_3^2 K^\exchF_3(X_3,X_3,X_3) \\
                     &= 8\cdot\frac{1}{4} - 4  + 2 = 0 \\
  \end{align*}
  and for all $n\ge4$ we clearly have $K^\exchF_n(Y_{i_1},Y_{i_2},\dots,Y_{i_n}) = 0$.
\end{proof}

\begin{Definition}
  Let $\exchm=(\alg{U},\tilde{\phi},\alg{J})$ be an exchangeability system
  for~$(\alg{A},\phi)$.
  We say that \emph{Marcinkiewicz' theorem holds} in \exch{} 
  if a selfadjoint random variable $X$,
  which for some fixed $m$ satisfies
  $K^\exchm_\pi(X_1,\dots,X_n)=0$ whenever one of the blocks of $\pi$ consists
  of more than $m$ copies of $X$'s, must be gaussian.

  We say that \emph{Cram\'er's theorem holds} in \exch{} if for any decomposition
  $X=X_1+X_2+\dots+X_n$ of a selfadjoint gaussian random variable $X$ into
  \exch-independent random variables $X_1,X_2,\dots,X_n$, the
  summands $X_j$ themselves must be gaussian.
\end{Definition}
\begin{Theorem}
  In an arbitrary exchangeability system
  Marcinkiewicz' theorem and Cram\'er's theorem imply the Skitovi\v{c}-Darmois theorem.
\end{Theorem}
\begin{proof}
  We follow the proof of \cite{KaganLinnikRao:1973:characterization}.
  Let $X_1,X_2,\dots,X_n$ be \exch-independent random variables
  and $a_i$, $b_i$ nonzero real numbers.
  Assume that the random variables $Y_1 = a_1 X_1 + a_2X_2 +\dots+a_nX_n$
  and $Y_2 = b_1 X_1 + b_2X_2 +\dots+b_nX_n$
  are \exch-independent.
  After rescaling the $X_i$'s we may assume that $a_i=1$ for all $i$.
  The assumed independence relations imply
  for every pair of real numbers $(\alpha,\beta)$ the following identities for
  the cumulants.
  \begin{align}
    K^\exchm_m(\alpha Y_1 + \beta Y_2) 
    &= \alpha^m K^\phi_m(Y_1) + \beta^m K^\phi_m(Y_2)\\
    &= \sum_{j=1}^n (\alpha^m + \beta^m b_j^m) K^\exchm_m(X_j)\\
    K^\exchm_m(\alpha Y_1 + \beta Y_2) 
    &= \sum_{j=1}^n (\alpha + \beta b_j)^m K^\exchm_m(X_j)\\
  \end{align}
  Let us first consider the case that the $b_j$ are pairwise different.
  Then we can differentiate the identity
  $$
  \sum_{j=1}^n (\alpha^m + \beta^m b_j^m) K^\exchm_m(X_j)
  =
  \sum_{j=1}^n (\alpha + \beta b_j)^m K^\exchm_m(X_j)
  $$
  $k$ times with respect to $\beta$ and evaluate at $\alpha=1$ and $\beta=0$ and obtain
  $$
  m(m-1)(m-2)\cdots(m-k+1) \sum b_j^k K^\exchm_m(X_j) = 0
  $$
  if $m>n$, this gives rise to a regular Vandermonde system and therefore
  $K^\exchm_m(X_j)$ must vanish.
  Marcinkiewicz' theorem then implies that the $X_j$ are gaussian.
  
  If some of the $b_j$'s are equal, we can group the $X_j$'s with equal coefficients
  together, 
  and the considerations above imply that the sum of the $X_j$'s in each group is
  gaussian. Then Cram\'er's theorem implies that the individual $X_j$'s are gaussian.
\end{proof}
As noted above, Marcinkiewicz' theorem and Cram\'er's theorems do not hold for free
independence, but they do hold e.g.\ for boolean independence,
see~\cite{SpeicherWoroudi:1997:boolean}.

\ssection{Quadratic forms. Lukacs' theorem}
The next result is known as \emph{Lukacs' theorem} in classical probability
\cite[\S{}III.6]{Feller:1971:introductionII}.
\begin{Proposition}
  Let~$X_1,X_2,\dots,X_n$ be a sequence of noncommutative \exch-i.i.d.\ random variables
  for which the singleton condition of 
  Definition~\ref{def:GeneralCumulants:Gaussian:SingletonCondition}
  holds. Then $X_i$ are gaussian if and only if
  their sample mean~$S_1=\sum X_k$ and sample variation
  $T=\sum (X_k-\frac{1}{n} S_1)^2 = \sum X_k^2- \frac{1}{n} S^2$
  are \exch-independent.
\end{Proposition}
\begin{proof}
  To prove necessity, choose any $n\times n$ orthogonal matrix $U=[U_{ij}]$
  with first row $U_{1j}=\frac{1}{\sqrt{n}}$.
  Then by Maxwell's theorem~\ref{thm:GeneralCumulants:Maxwell} the random variables
  $X_i'=\sum U_{ij} X_j$ are also \exch-independent.
  Consequently $S_1=\sqrt{n} X_1'$ and
  $$
  T=\sum X_k^2-\frac{1}{n} S_1^2=\sum X_k'{}^2-X_1'{}^2
   =\sum_2^n X_j'{}^2
  $$
  are \exch-independent.

  In order to prove sufficiency of the condition,
  we show that the presence of a block of length at least three in a partition
  implies that the corresponding cumulant vanishes.
  It is enough to consider the full cumulants, the argument for partitioned
  cumulants is entirely similar.
  Thus assume $m\geq3$, then the \exch-independence of $S_1$ and $T=S_2-\frac{1}{n}S_1^2$ implies
  \begin{align*}
    0 &= K^\exchm_{m-1}(S_1,S_1,\dots,S_1,T)\\
      &= K^\exchm_{m-1}(S_1,S_1,\dots,S_1,S_2) - \frac{1}{n}\,K^\exchm_{m-1}(S_1,S_1,\dots,S_1,S_1^2)
  \end{align*}
  By the product formula
  (Proposition~\externalref{I}{prop:GeneralCumulants:Productformula})
  we have 
  $$
  K^\exchm_{m-1}(Y,Y,\dots,Y,Y^2)
  = K^\exchm_m(Y,Y,\dots,Y)
    + \sum_{\substack{\pi\vee\pi_0=\hat1_{m}\\ \pi<\hat1_m}} K^\exchm_\pi(Y,Y,\dots,Y)
  $$
  for any random variable $Y$, where $\pi_0=
  \begin{picture}(80,8.4)(1,0)
    \put(10,0){\line(0,1){8.4}}
    \put(20,0){\line(0,1){8.4}}
    \put(30,0){\line(0,1){8.4}}
    \put(40,0){$\cdots$}
    \put(60,0){\line(0,1){8.4}}
    \put(70,0){\line(0,1){8.4}}
    \put(80,0){\line(0,1){8.4}}
    \put(10,8.4){\line(1,0){0}}
    \put(20,8.4){\line(1,0){0}}
    \put(30,8.4){\line(1,0){0}}
    \put(40,8.4){\line(1,0){0}}
    \put(50,8.4){\line(1,0){0}}
    \put(60,8.4){\line(1,0){0}}
    \put(70,8.4){\line(1,0){10}}
  \end{picture}
  $.
  Therefore the first term is
  \begin{align*}
  K^\exchm_{m-1}(S_1,S_1,\dots,S_1,S_2)
  &= n  K^\exchm_{m-1}(X,X,\dots,X,X^2)\\
  &= n
     \biggl(
       K^\exchm_m(X,X,\dots,X)
       + \sum_{\substack{\pi\vee\pi_0=\hat1_{m}\\ \pi<\hat1_m}} K^\exchm_\pi(X,X,\dots,X)
     \biggr)
  \end{align*}
  while the second term is
  \begin{align*}
    \frac{1}{n}\,K^\exchm_{m-1}(S_1,S_1,\dots,S_1,S_1^2) 
    &= \frac{1}{n}
       \biggl(
         K^\exchm_m(S_1,S_1,\dots,S_1)
         + 
         \sum_{\substack{\pi\vee\pi_0=\hat1_{m}\\ \pi<\hat1_m}} K^\exchm_\pi(S_1,S_1,\dots,S_1)
       \biggr)\\
    &= K^\exchm_m(X,X,\dots,X)
       + n \sum_{\substack{\pi\vee\pi_0=\hat1_{m}\\ \pi<\hat1_m}} K^\exchm_\pi(X,X,\dots,X)
  \end{align*}
  because each partition $\pi$ in the sum has exactly two blocks;
  The difference of the two terms is $(n-1)\,K^\exchm_m(X,X,\dots,X)$ and vanishes.
\end{proof}

\ddate{19.09.2002}
Lukacs' theorem can be generalized to more general quadratic forms as follows.
For different proofs in the free case
see~\cite{HiwatashiNagisaYoshida:1999:characterizations,HiwatashiKurodaNagisaYoshida:1999:freeanalogue}.
We are grateful to H.~Yoshida for bringing the latter to our attention.
\begin{Proposition}[{\cite[Prop.~2.2]{HiwatashiNagisaYoshida:1999:characterizations}}]
  Let $X_i$ be a sequence of \exch-i.i.d.\ Gaussian random variables in the sense
  of Definition~\ref{def:GeneralCumulants:Gaussian:SameDistribution}
  and
  let $A\in M_n(\IR)$, $b\in\IR^n$ such that
  \begin{equation}
    \label{eq:GeneralCumulants:Gaussian:Ab=0andbtA=0}
    A b = 0 \qquad\qquad b^t A=0
  \end{equation}
  Then the linear form $L=\sum b_i X_i$ and the quadratic form
  $Q=\sum a_{ij} X_i X_j$ are \exch-independent.
\end{Proposition}
\begin{proof}
  Without loss of generality we may assume that $b$ is a unit vector,
  i.e., $b^t b=1$.
  In that case we can extend $b$ to an orthonormal basis of $\IR^n$,
  denoted $\{b_1=b,b_2,\dots,b_n\}$ where $b_i$ has components $(b_{ij})_{j=1,\dots,n}$.
  We can express $A$ in this basis as
  $$
  A = \sum \alpha_{ij} b_i b_j^t
  $$
  with $\alpha_{ij} = b_i^t A b_j$.
  By Maxwell's Theorem~\ref{thm:GeneralCumulants:Maxwell}
  the sequence $Y_i = \sum b_{ij}X_j$ has the same distribution as $X_i$.
  Our assumption~\eqref{eq:GeneralCumulants:Gaussian:Ab=0andbtA=0}
  implies that $\alpha_{1j}=0$ and $\alpha_{j1}=0$ for all $j$
  and we can rewrite $L$ and $Q$ as
  $$
  L = Y_1
  \qquad\qquad
  Q = \sum_{i,j=2}^n \alpha_{ij} Y_i Y_j
  $$
  which are clearly \exch-independent.
\end{proof}

\begin{Proposition}[{\cite[Thm.~2.3]{HiwatashiNagisaYoshida:1999:characterizations}}]
  Let~$(X_i)$ be an \exch-i.i.d.\ centered sequence satisfying the singleton condition
  and let $A=[a_{ij}]\in M_n(\IR)$,
  $b\in \IR^n$ be such that
  $$
  A b=b^t A=0
  \qquad\qquad
  \forall m\in\IN:\sum_i b_i^m a_{ii} \ne 0
  .
  $$
  If $L=\sum b_i X_i$ and $Q=\sum a_{ij} X_i X_j$ are \exch-independent,
  then $X_i$ are Gaussian.
\end{Proposition}
\begin{proof}
  The singleton condition implies that cumulants with singleton blocks vanish.
  Therefore it suffices to show that cumulants with a block of length greater than or
  equal to three vanish.
  We proceed by induction.
  Consider for $m\geq2$ the cumulant
  $$
  K^\exchm_m(L,L,\dots,L,Q) = 0.
  $$
  We can expand it with the help of the product
  formula (Proposition~\externalref{I}{prop:GeneralCumulants:Productformula}):
  \begin{align*}
    K^\exchm_m&(L,L,\dots,L,Q)\\
    &= \sum_{h:[m+1]\to[n]}
        b_{h(1)}
        b_{h(2)}
        \cdots
        b_{h(m-1)}
        a_{h(m)h(m+1)}
        K^\exchm_m(X_{h(1)},\dots,X_{h(m-1)},X_{h(m)}X_{h(m+1)}) \\
    &= \sum_{\pi\in\Pi_{m+1}}
       \sum_{\substack{h:[m+1]\to[n] \\ \ker h=\pi}}
        b_{h(1)}
        b_{h(2)}
        \cdots
        b_{h(m-1)}
        a_{h(m)h(m+1)}
        K^\exchm_m(X_{\pi(1)},\dots,X_{\pi(m-1)},X_{\pi(m)}X_{\pi(m+1)}) \\
    \intertext{
      Note that each contributing partition~$\pi$ has at most two blocks
      and by the product formula
      $$
      K^\exchm_m(X_{\pi(1)},\dots,X_{\pi(m-1)},X_{\pi(m)}X_{\pi(m+1)})
      = \sum_{\substack{\rho\vee \pi_0=\hat{1}_{m+1} \\ \rho\leq\pi }}
         K^\exchm_\rho(X,X,\dots,X)
      $$
      with $\pi_0 =
      \begin{picture}(80,8.4)(1,0)
        \put(10,0){\line(0,1){8.4}}
        \put(20,0){\line(0,1){8.4}}
        \put(30,0){\line(0,1){8.4}}
        \put(40,0){$\cdots$}
        \put(60,0){\line(0,1){8.4}}
        \put(70,0){\line(0,1){8.4}}
        \put(80,0){\line(0,1){8.4}}
        \put(10,8.4){\line(1,0){0}}
        \put(20,8.4){\line(1,0){0}}
        \put(30,8.4){\line(1,0){0}}
        \put(40,8.4){\line(1,0){0}}
        \put(50,8.4){\line(1,0){0}}
        \put(60,8.4){\line(1,0){0}}
        \put(70,8.4){\line(1,0){10}}
      \end{picture}
      $,
      and where $X=X_1$ has the same distribution as all the $X_i$.
      Using this we can continue
    }
    &= \sum_{\substack{\pi,\rho\in\Pi_{m+1}\\\rho\vee \pi_0=\hat{1}_{m+1} \\ \rho\leq\pi} }
       \sum_{\ker h=\pi}
        b_{h(1)}
        b_{h(2)}
        \cdots
        b_{h(m-1)}
        a_{h(m)h(m+1)}
        K^\exchm_\rho(X,X,\dots,X)\\
    &= \sum_{\rho\vee \pi_0=\hat{1}_{m+1} }
       \sum_{\ker h\geq\rho}
        b_{h(1)}
        b_{h(2)}
        \cdots
        b_{h(m-1)}
        a_{h(m)h(m+1)}
        K^\exchm_\rho(X,X,\dots,X) \\
    &= \sum_i b_i^{m-1} a_{ii} K^\exchm_{m+1}(X,X,\dots,X) \\
    &\phantom{===}   +
       \sum_{\substack{\rho\vee \pi_0=\hat{1}_{m+1} \\ \rho < \hat{1}_{m+1} }}
       \sum_{\ker h\geq\rho}
        b_{h(1)}
        b_{h(2)}
        \cdots
        b_{h(m-1)}
        a_{h(m)h(m+1)}
        K^\exchm_\rho(X,X,\dots,X)
  \end{align*}
  We will now apply induction to show that all but the first term of the last summation
  vanish, and together with the assumption $\sum b_i^{m-1} a_{ii}\ne0$ this will
  imply that $K^\exchm_{m+1}(X)=0$.

  For $m=2$ we have
  \begin{align*}
    0 = K^\exchm_2(L,Q)
    &= \sum b_i a_{ii} K^\exchm_3(X,X,X) \\ &\phantom{===}
       +
       \sum_{\rho\in
         \{\!\begin{picture}(15,4.2)(1,0)
             \put(5,0){\line(0,1){4.2}}
             \put(10,0){\line(0,1){4.2}}
             \put(15,0){\line(0,1){4.2}}
             \put(5,4.2){\line(1,0){5}}
             \put(15,4.2){\line(1,0){0}}
           \end{picture}\,,\!
           \begin{picture}(15,7.2)(1,0)
             \put(5,0){\line(0,1){7.2}}
             \put(10,0){\line(0,1){4.2}}
             \put(15,0){\line(0,1){7.2}}
             \put(5,7.2){\line(1,0){10}}
             \put(10,4.2){\line(1,0){0}}
           \end{picture}\,
         \}}
        \sum_{\ker h\geq\rho}
         b_{h(1)}
         a_{h(2)h(3)}
         K^\exchm_\rho(X,X,X)
  \end{align*}
  and because of the singleton condition all but the first summand vanish,
  showing that $K^\exchm_3(X,X,X)=0$.
  
  For $m=3$ we have
  \begin{align*}
    0 = K^\exchm_3(L,L,Q)
    &= \sum b_i^2 a_{ii} K^\exchm_4(X,X,X,X) \\ &\phantom{===}
       +
       \sum_{\rho\in
         \{\!
           \begin{picture}(20,4.2)(1,0)
             \put(5,0){\line(0,1){4.2}}
             \put(10,0){\line(0,1){4.2}}
             \put(15,0){\line(0,1){4.2}}
             \put(20,0){\line(0,1){4.2}}
             \put(5,4.2){\line(1,0){10}}
             \put(20,4.2){\line(1,0){0}}
           \end{picture}
           \,,\!
           \begin{picture}(20,7.2)(1,0)
             \put(5,0){\line(0,1){7.2}}
             \put(10,0){\line(0,1){7.2}}
             \put(15,0){\line(0,1){4.2}}
             \put(20,0){\line(0,1){7.2}}
             \put(5,7.2){\line(1,0){15}}
             \put(15,4.2){\line(1,0){0}}
           \end{picture}
           \,,\!
           \begin{picture}(20,7.2)(1,0)
             \put(5,0){\line(0,1){4.2}}
             \put(10,0){\line(0,1){7.2}}
             \put(15,0){\line(0,1){4.2}}
             \put(20,0){\line(0,1){7.2}}
             \put(5,4.2){\line(1,0){10}}
             \put(10,7.2){\line(1,0){10}}
           \end{picture}
           \,,\!
           \begin{picture}(20,7.2)(1,0)
             \put(5,0){\line(0,1){7.2}}
             \put(10,0){\line(0,1){4.2}}
             \put(15,0){\line(0,1){4.2}}
             \put(20,0){\line(0,1){7.2}}
             \put(10,4.2){\line(1,0){5}}
             \put(5,7.2){\line(1,0){15}}
           \end{picture}\,
         \}}
        \sum_{\ker h\geq\rho}
         b_{h(1)}
         b_{h(2)}
         a_{h(3)h(4)}
         K^\exchm_\rho(X,X,X,X)
  \end{align*}
  For 
  $\rho\in
  \{\!
  \begin{picture}(20,4.2)(1,0)
    \put(5,0){\line(0,1){4.2}}
    \put(10,0){\line(0,1){4.2}}
    \put(15,0){\line(0,1){4.2}}
    \put(20,0){\line(0,1){4.2}}
    \put(5,4.2){\line(1,0){10}}
    \put(20,4.2){\line(1,0){0}}
  \end{picture}
  \,,\!
  \begin{picture}(20,7.2)(1,0)
    \put(5,0){\line(0,1){7.2}}
    \put(10,0){\line(0,1){7.2}}
    \put(15,0){\line(0,1){4.2}}
    \put(20,0){\line(0,1){7.2}}
    \put(5,7.2){\line(1,0){15}}
    \put(15,4.2){\line(1,0){0}}
  \end{picture}
  \,\}$,
  the term vanishes by induction hypothesis (and because of the singleton condition).
  For
  $\rho\in\{\!
  \begin{picture}(20,7.2)(1,0)
    \put(5,0){\line(0,1){4.2}}
    \put(10,0){\line(0,1){7.2}}
    \put(15,0){\line(0,1){4.2}}
    \put(20,0){\line(0,1){7.2}}
    \put(5,4.2){\line(1,0){10}}
    \put(10,7.2){\line(1,0){10}}
  \end{picture}
  \,,\!
  \begin{picture}(20,7.2)(1,0)
    \put(5,0){\line(0,1){7.2}}
    \put(10,0){\line(0,1){4.2}}
    \put(15,0){\line(0,1){4.2}}
    \put(20,0){\line(0,1){7.2}}
    \put(10,4.2){\line(1,0){5}}
    \put(5,7.2){\line(1,0){15}}
  \end{picture}\,
  \}$
  the coefficient of $K^\exchm_\rho$ is $\sum b_i b_j a_{ij}$
  which vanishes by assumption on $A$ and $b$.
  Now for $m\geq4$ any $\rho< \hat{1}_{m+1}$ satisfying $\rho\vee\pi_0=\hat{1}_{m+1}$
  has exactly two blocks and one of the blocks has cardinality at least three,
  and the induction hypothesis implies that the cumulant vanishes.
\end{proof}

\begin{Proposition}[{\cite[Prop.~2.2]{HiwatashiKurodaNagisaYoshida:1999:freeanalogue}}]
  Let $X_i$ be \exch-i.i.d.\ copies of the Gaussian random variable~$X$.
  Then the quadratic form
  $Q=\sum a_{ij} X_i X_j$ with $A=[a_{ij}]$ symmetric has cumulants
  $$
  K^\exchm_n(Q) = \tr(A^n)\, K^\exchm_n(X^2)
  $$
\end{Proposition}
\begin{proof}
  Indeed~$A$ can be diagonalized to $A=U^t\Lambda U$ with~$U$ orthogonal
  and by Proposition~\ref{thm:GeneralCumulants:Maxwell}
  the random variables~$Y_i=\sum U_{ij} X_j$ have the same distribution as $X_i$
  and the cumulants of $\sum \lambda_i Y_i^2$ are
  $$
  K^\exchm_n( \sum \lambda_i Y_i^2 ) = \sum \lambda_i^n K^\exchm_n(Y_i^2)
  $$
\end{proof}

\begin{Remark}
  The joint cumulants of arbitrary quadratic forms in Gaussian random variables
  are computed as follows.
  Let~$X_i$ be an \exch-i.i.d.\ sequence of a Gaussian random variable~$X$ and
  let $Q_k = \sum a_{ij}(k) X_i X_j$ be quadratic forms
  where $A_k = [a_{ij}(k)]$ are not necessarily symmetric matrices.
  Then
  \begin{align*}
    K^\exchm_m&(Q_1,Q_2,\cdots,Q_m)\\
    &= \sum_{\pi\in\Pi_{2m}}
        \sum_{\ker h=\pi} 
         a_{h(1)h(2)}(1)\cdots a_{h(2m-1)h(2m)}(m)
         K^\exchm_m(X_{\pi(1)}X_{\pi(2)},\dots,X_{\pi(2m-1)}X_{\pi(2m)})\\
    &= \sum_{\substack{\rho\vee\pi_0 = \hat{1}_{2m}\\
                        \rho\in\Pi_{2m}^{(2)}}}
        \sum_{\ker h\geq \rho}
         a_{h(1)h(2)}(1)\cdots a_{h(2m-1)h(2m)}(m)
         K^\exchm_\rho(X)
  \end{align*}
  where~$
  \pi_0=
  \begin{picture}(45,4.2)(1,0)
    \put(5,0){\line(0,1){4.2}}
    \put(10,0){\line(0,1){4.2}}
    \put(15,0){\line(0,1){4.2}}
    \put(20,0){\line(0,1){4.2}}
    \put(25,0){\!$\cdots$}
    \put(40,0){\line(0,1){4.2}}
    \put(45,0){\line(0,1){4.2}}
    \put(5,4.2){\line(1,0){5}}
    \put(15,4.2){\line(1,0){5}}
    \put(25,4.2){\line(1,0){0}}
    \put(30,4.2){\line(1,0){0}}
    \put(35,4.2){\line(1,0){0}}
    \put(40,4.2){\line(1,0){5}}
  \end{picture}
  $.
  For free Gaussians (that is, free semicircular random variables)
  there is only one contributing partition, namely
  $\rho_0=
  \begin{picture}(55,7.2)(1,0)
    \put(5,0){\line(0,1){7.2}}
    \put(10,0){\line(0,1){4.2}}
    \put(15,0){\line(0,1){4.2}}
    \put(20,0){\line(0,1){4.2}}
    \put(25,0){\line(0,1){4.2}}
    \put(30,0){\!$\cdots$}
    \put(45,0){\line(0,1){4.2}}
    \put(50,0){\line(0,1){4.2}}
    \put(55,0){\line(0,1){7.2}}
    \put(10,4.2){\line(1,0){5}}
    \put(20,4.2){\line(1,0){5}}
    \put(40,4.2){\line(1,0){0}}
    \put(45,4.2){\line(1,0){5}}
    \put(5,7.2){\line(1,0){50}}
\end{picture}
  $,
  and in this case
  $$
  K^\exchF_m(Q_1,Q_2,\cdots,Q_m) = \tr(A_1 A_2\cdots A_m)
  .
  $$
  For classical Gaussian random variables,
  every pair partition~$\rho$ with~$\rho\vee\rho_0=\hat{1}_{2m}$
  can be obtained from~$\rho_0$ by permuting and flipping the pairs 
  $(3,4),(5,6),\dots,(2m-1,2m)$. A permutation of pairs corresponds to
  a permutation of the matrices $A_k$ and a flip corresponds to replacing
  $A_k$ by its transpose $A_k^t$. 
  There are $(n-1)!\cdot 2^{n-1}$ ways to do this
  and we get
  $$
  \kappa_m(Q_1,Q_2,\cdots,Q_m)
  = \sum_{\sigma\in\SG_{\{2,3,\dots,m\}}}
    \sum_{\eps_2,\eps_3,\dots,\eps_m\in\{1,t\}}
     \tr(A_1 
         A_{\sigma(2)}^{\eps_2}
         A_{\sigma(3)}^{\eps_3}
         \cdots
         A_{\sigma(m)}^{\eps_m})
  $$
  In the general case, these summands are weighted with the corresponding
  $K^\exchm_\rho(X)$.
\end{Remark}

\begin{Proposition}[{\cite[Prop.~2.3]{HiwatashiKurodaNagisaYoshida:1999:freeanalogue}}]
  Let~$X_i$ be an \exch-i.i.d.\ Gaussian sequence and
  $A$, $B\in M_n(\IR)$ symmetric matrices.
  Then the quadratic forms
  $$
  Q = \sum a_{ij} X_i X_j
  \qquad\qquad
  Q' = \sum b_{ij} X_i X_j
  $$
  are \exch-independent if and only if $AB=0$.
\end{Proposition}
\begin{proof}
  Assume that~$Q$ and~$Q'$ are \exch-independent.
  Then we can write cumulants in two ways:
  $$
  K^\exchm_n(s Q + t Q')
  = \tr ( (sA+tB)^n ) K^\exchm_n(X^2)
  = (s^n \tr(A^n) + t^n \tr(B^n) ) K^\exchm_n(X^2)
  $$
  for~$n=4$ this implies
  $$
  \tr (( AB+BA )^2 ) + 2\tr( BA^2B) = 0
  $$
  and therefore $\tr ((AB)^t(AB))=0$, i.e., $AB=0$.

  For the converse, we will prove that more generally,
  the symmetric quadratic forms $Q_k=\sum a_{ij}(k) X_i X_j$ are \exch-independent
  if the matrices~$A_k=[a_{ij}(k)]$ satisfy $A_i A_j=0$ for $i\ne j$.
  Indeed, by the product formula
  (Proposition~\externalref{I}{prop:GeneralCumulants:Productformula})
  we have
  \begin{align*}
    K_m&\exchm(Q_{k_1},\dots,Q_{k_m})\\
    &=  \sum_{h:[2m]\to [n]}
         a_{h(1)h(2)}(k_1)
         \cdots
         a_{h(2m-1)h(2m)}(k_m)
         \,
         K^\exchm_m(X_{h(1)} X_{h(2)},\dots, X_{h(2m-1)}X_{h(2m)})\\
    &= \sum_{\pi\in\Pi_{2m}}
        \sum_{\ker h=\pi}
         a_{h(1)h(2)}(k_1)
         \cdots
         a_{h(2m-1)h(2m)}(k_m)
         \,
         K^\exchm_m(X_{\pi(1)} X_{\pi(2)},\dots, X_{\pi(2m-1)}X_{\pi(2m)})\\
    &= \sum_{\pi\in\Pi_{2m}}
        \sum_{\ker h=\pi}
         a_{h(1)h(2)}(k_1)
         \cdots
         a_{h(2m-1)h(2m)}(k_m)
         \sum_{\substack{\rho\leq\pi \\ \rho\vee \rho_0 = \hat{1}_{2m}}}
          K^\exchm_\rho(X) \\
    &= \sum_{\substack{\rho\in\Pi_{2m}^{(2)} \\ \rho\vee \rho_0 = \hat{1}_{2m}}}
        \sum_{\ker h\geq\rho}
         a_{h(1)h(2)}(k_1)
         \cdots
         a_{h(2m-1)h(2m)}(k_m)
         K^\exchm_\rho(X)
  \end{align*}
  where $\rho_0=
  \begin{picture}(45,4.2)(1,0)
    \put(5,0){\line(0,1){4.2}}
    \put(10,0){\line(0,1){4.2}}
    \put(15,0){\line(0,1){4.2}}
    \put(20,0){\line(0,1){4.2}}
    \put(25,0){\!$\cdots$}
    \put(40,0){\line(0,1){4.2}}
    \put(45,0){\line(0,1){4.2}}
    \put(5,4.2){\line(1,0){5}}
    \put(15,4.2){\line(1,0){5}}
    \put(35,4.2){\line(1,0){0}}
    \put(40,4.2){\line(1,0){5}}
  \end{picture}
  $.
  Each~$\rho$ is a pair partition in which each block connects
  two different blocks of $\rho_0$ in such a way that the resulting
  graph is connected. If we number the blocks of $\rho_0$ from $1$ to $m$,
  we can define a cycle $\sigma\in\SG_m$ starting at block $1$, choosing an arc of
  $\rho$ which connects it so some block $\sigma(1)$, choosing the other arc
  starting in block $\sigma(1)$ etc.n
  Using the symmetry of the matrices $A_k$ we can rewrite
  \begin{align*}
  \sum_{\ker h\geq\rho}
   a_{h(1)h(2)}(k_1)
   \cdots
   a_{h(2m-1)h(2m)}(k_m)
   &= \sum_{i_1,\dots,i_m}
       a_{i_{\rho(1)} i_{\rho(2)}}(k_1)
       \cdots
       a_{i_{\rho(2m-1)} i_{\rho(2m)}}(k_m)\\
   &= \tr(A_{k_{\sigma(1)}} A_{k_{\sigma(2)}} \cdots A_{k_{\sigma(m)}})
   .
  \end{align*}
  If the $k_i$ are not all equal,
  the matrix product vanishes by assumption.
\end{proof}

\begin{Proposition}[{\cite[Thm.~3.2]{HiwatashiKurodaNagisaYoshida:1999:freeanalogue}}]
  Let $X_i$ be a sequence of \exch-i.i.d.\ copies of the Gaussian random variable $X$.
  Then the distribution of
  $$
  Y = \sum_{i=1}^n (X_i+a_i)^2
  $$
  only depends on the number $\sum a_i^2$.
\end{Proposition}
\begin{proof}
  We show that the cumulants only depend on $\sum a_i^2$.
  \begin{align*}
    K^\exchm_m(Y)
    &= \sum_i K^\exchm_m( (X+a_i)^2) \\
    &= \sum_i \sum_{\pi\vee \pi_0 = \hat{1}_{2m}} K^\exchm_\pi(X+a_i)\\
\intertext{
  where
  $\pi_0 = 
  \begin{picture}(45,4.2)(1,0)
    \put(5,0){\line(0,1){4.2}}
    \put(10,0){\line(0,1){4.2}}
    \put(15,0){\line(0,1){4.2}}
    \put(20,0){\line(0,1){4.2}}
    \put(25,0){\!$\cdots$}
    \put(40,0){\line(0,1){4.2}}
    \put(45,0){\line(0,1){4.2}}
    \put(5,4.2){\line(1,0){5}}
    \put(15,4.2){\line(1,0){5}}
    \put(35,4.2){\line(1,0){0}}
    \put(40,4.2){\line(1,0){5}}
  \end{picture}
  $.
  Since~$X$ is Gaussian, the blocks of each contributing~$\pi$ have size at most~$2$.
  Together with the condition~$\pi\vee\pi_0=\hat{1}_{2m}$ this implies that $\pi$ is
  either a pair partition or~$\pi$ has $2$ singletons.
}
    &= \sum_i
        \biggl(
          \sum_{\substack{ \pi\vee \pi_0 = \hat{1}_{2m} \\ \pi\in\Pi_{2m}^{(2)}}}
           K^\exchm_\pi(X+a_i)
          +
          \sum_{\substack{ \pi\vee \pi_0 = \hat{1}_{2m} \\ \pi\not\in\Pi_{2m}^{(2)}}}
           K^\exchm_\pi(X+a_i)
        \biggr)
\intertext{
  In the first term there are only cumulants of order~$2$ which are invariant
  under translations, therefore we can forget the $a_i$.
  In the second term, there are exactly~$2$ singletons 
  and by Lemma~\externalref{I}{lem:GeneralCumulants:RemoveIdentity}
  we have for such partitions~$\pi$ that
  $$
  K^\exchm_\pi(X+a_i) = a_i^2 K^\exchm_{\tilde{\pi}}(X)
  $$
  where $\tilde{\pi}\in\Pi_{2m-2}^{(2)}$ is the partition obtained from
  $\pi$ by removing the two singletons.
}
    &= \sum_i
        \biggl(
          K^\exchm_m(X^2)
          +
          \sum_{\substack{ \pi\vee \pi_0 = \hat{1}_{2m} \\ \pi\not\in\Pi_{2m}^{(2)}}}
           a_i^2
           K^\exchm_{\tilde{\pi}}(X)
        \biggr)\\
    &= n K^\exchm_m(X^2)
       +
       \left(
         \sum a_i^2
       \right)
       \sum_{\substack{ \pi\vee \pi_0 = \hat{1}_{2m} \\ \pi\not\in\Pi_{2m}^{(2)}}}
        K^\exchm_{\tilde{\pi}}(X)
  \end{align*}
\end{proof}
A converse of this theorem holds in classical and free probability,
but we were not able to find a generalization.

\emph{Acknowledgements}.
We are grateful to Hiroaki Yoshida for bringing to our attention the papers
\cite{HiwatashiNagisaYoshida:1999:characterizations,HiwatashiKurodaNagisaYoshida:1999:freeanalogue}
and we particularly thank Philippe Biane for a careful reading of the manuscript.
